\documentclass{elsart}
 \usepackage{graphics}
 \usepackage{graphicx}
 \usepackage{epsfig}
\usepackage{amssymb}
\newtheorem{lemma}{Lemma}[section]
\newtheorem{theorem}{Theorem}[section]

\begin{document}
\begin{frontmatter}
\title{Classification of the solutions of a mixed nonlinear Schr\"odinger system}

\author{Riadh Chteoui}
\ead{rchteoui@yahoo.fr}
\address{Laboratory of Algebra, Number Theory and Nonlinear Analysis LR18ES15, Department of Mathematics, Faculty of Science, University of Monastir, Tunisia.}
\author{Mohamed Lakdar Ben Mohamed}
\ead{lakdar.mohamed@fsm.rnu.tn; mbenmohamed@kku.edu.sa}
\address{Laboratory of Algebra, Number Theory and Nonlinear Analysis LR18ES15, Department of Mathematics, Faculty of Science, University of Monastir, Tunisia.\\
and\\
Department of Computer Sciences, Higher Institute of Applied Sciences and Technology, University of Sousse, City Ibn Khaldoun 4003, Sousse, Tunisia.\\
and\\
Department of Mathematics, Faculty of Sciences and Arts, University of King Khaled, Mahayel Branch, Saudi Arabia.}
\author{Abdulrahman F. Aljohani}
\ead{a.f.aljohani@ut.edu.sa}
\address{Department of Mathematics, Faculty of Science, University of Tabuk, Saudi Arabia.}
\author{Anouar Ben Mabrouk\corauthref{cor1}}
\ead{anouar.benmabrouk@fsm.rnu.tn; amabrouk@ut.edu.sa}
\corauth[cor1]{Corresponding author}
\address{Laboratory of Algebra, Number Theory and Nonlinear Analysis LR18ES15, Department of Mathematics, Faculty of Sciences, University of Monastir, Tunisia.\\
and\\
Department of Mathematics, Higher Institute of Applied Mathematics and Informatics, University of Kairouan, Tunisia.\\
and\\	
Department of Mathematics, Faculty of Science, University of Tabuk, Saudi Arabia.
}
\begin{abstract}
In this paper, we study a couple of NLS equations characterized by mixed cubic and superlinear power laws. Classification of the solutions as well as existence and uniqueness of the steady state solutions have been investigated. 
\end{abstract}
\end{frontmatter}
\section{Introduction}
The paper is devoted to the study of some nonlinear systems of PDEs of the form
\begin{equation}\label{NLCP}
\begin{cases}
\mathcal{L}_1(u)+f_1(u,v)=0,\\
\mathcal{L}_2(v)+f_2(u,v)=0
\end{cases}
\end{equation}
where $u=u(x,t)$, $v=v(x,t)$ on a domain $\Omega\times(t_0,+\infty)$, $\Omega\subset\mathbb{R}^N$, $t_0\in\mathbb{R}$ with suitable initial and boundary conditions. The operators $\mathcal{L}_1$ and $\mathcal{L}_2$ are linear Schr\"odinger-type operators of the form
$$
\mathcal{L}_i(u(x,t))=iu_t+\sigma_iu_{xx},\;i=1,2
$$ 
leading to a nonlinear Schr\"odinger system. $u_{xx}$ is the second order partial derivative relatively to $x$, which plays the role of the Laplacian, $u_t$ is the first order partial derivative in time, $\sigma_i$, $i=1,2$ are constant real positive parameters. $f_1$ and $f_2$ are nonlinear continuous functions of two variables. 

Remark that whenever the $\mathcal{L}_i$'s are the classical Schr\"odinger operators
$$
\mathcal{L}_1(u)=\mathcal{L}_2(u)=iu_t+u_{xx}
$$
and $u=v=\varphi$, $f_1=f_2=f$, we come back to the original Schr\"odinger equation
\begin{equation}\label{NLS-phi}
i\varphi_t+\Delta \varphi=f(\varphi).
\end{equation}
On the other hand, this last equation itself may lead to a system of PDEs of real valued functions satisfying a Heat system. Indeed, let $\varphi=u+iv=(1+i)u$, we get a system of coupled Heat ones
\begin{equation}\label{NLCP1}
\begin{cases}
u_t-\Delta u+f_1(u,v)=0.\\
v_t+\Delta v+f_2(u,v)=0,
\end{cases}
\end{equation}
where $f_1$ and $f_2$ are issued from the real and imagiary parts of $f(\varphi)$. 

Return again to the system (\ref{NLCP}) in the simple case $\sigma_i=1$ and denote $\varphi=u+iv$ and $f(\varphi)=f_1(u,v)+if_2(u,v)$, we come back to the classical NLS equation where $u$ and $v$ may be seen as real and imaginary parts of the solution $\varphi$ of an equation of the form (\ref{NLS-phi}).

As it is related to many physical/natural phenomena such as plasma, optics, condensed matter physics, etc, nonlinear Schr\"odinger system has attracted the interest of researchers in different fields such as pure and applied mathematics, pure and applied physics, quantum mechanics, mathematical physics, and it continues to attract resaerchers nowadays with the discovery of nano-physics, fractal domains, planets understanding, ... etc. For instence, in hydrodynamics, the NLS system may be a good model to describe the propagation of packets of waves according to some directions where a phenomenon of overlopping group velocity projection may occur \cite{A.K.Dhar}. In optics also, the propagation of short pulses has been investigated via a system of NLS equations \cite{Menyuk}. 

In single NLS equation, studies have been well developed from both theoretical and numerical aspects. The most recent mixed model is developed in \cite{BenmabroukAyadi-1}, \cite{BenmabroukAyadi-2}, \cite{BenmabroukBenMohamed-Nodal}, \cite{BenmabroukBenMohamed-Classification}, \cite{BenmabroukBenMohamed-NonRadial}, \cite{BenmabroukBenMohamed-FJAM}, \cite{BenmabroukBenMohamedOmrani}, \cite{ChteouiBenmabroukOunaies}, \cite{Kenig-Merle}, \cite{Keraani}, \cite{Martel-Merle}, \cite{Merle}, \cite{Yadava}, \cite{Yanagida},  ... with a general model
$$
f(u)=|u|^{p-1}u+\lambda|u|^{q-1}u.
$$
which coincides with $f(u,u)$ in the model (\ref{f-model}). An extending interesting model will be a mixed model of nonlinearities 
\begin{equation}\label{f-model}
f(u,v)=(|u|^{p-1}+\lambda|v|^{q-1})u,\;\;p,q,\lambda\in\mathbb{R}.
\end{equation}
In the literature, few works are done on such a model and are focusing on the mixed cubic-cubic ($p=q=3$) one 
$$
f(u,v)=(|u|^{2}+\lambda|v|^{2})u.
$$
For example, in \cite{Hai-Qiang.Zhang}, (2+1)-dimensional coupled NLS equations have been studied based on symbolic computation and Hirota method via the cubic-cubic nonlinear system
\begin{equation}\label{0.1}
\left\{\begin{array}{lll}
iu_t+u_{xx}+\sigma(|u|^2+\alpha|v|^2)u=0,\\
iv_t+v_{xx}+\sigma(|v|^2+\alpha|u|^2)v=0,
\end{array}\right.
\end{equation}
where $\alpha$ and $\sigma$ are real parameters. The same system has been also studied by many authors such as \cite{Aitchison}, \cite{Chakravarty}, \cite{Hioe}, \cite{Kanna}, \cite{Mollenauer}, \cite{Shalaby}. In \cite{Yan.Zhida}, the following more general system has been investigated for the asymptotic time behavior of the solutions,
\begin{equation}\label{2.2}
\left\{\begin{array}{lll}
i\partial_tu+\alpha\partial_x^2u+A(|u|^2+|v|^2)u=0,\\
i\partial_tv+\alpha\partial_x^2v+A(|u|^2+|v|^2)v=0\\
\end{array}\right.
\end{equation}
See also \cite{Gupta}. In \cite{K.Chaib} the following $p$-Laplacian stationary system has been discussed for necessary and sufficient conditions for the existence of the solutions
\begin{equation}\label{5.1}
\left\{
\begin{array}{c}
-\Delta_pu=\mu_1\Gamma_1(x,u,v) ~~\hbox{in} ~~\mathbb{R}^N\\
-\Delta_pv=\mu_2\Gamma_2(x,u,v)~~ \hbox{in}~~ \mathbb{R}^N,
\end{array}
\right.
\end{equation}
where $\Delta_pu=div(|\nabla u|^{p-2}\nabla u)$ is the $p$-Laplacian $(1<p<N)$. In \cite{Benci-Fortunato}, solutions of the type
\begin{equation}\label{6.2}
u(x,t)=u(x)e^{iwt},
\end{equation}
where $u$ is a real function known as standing wave or steady state and $w\in\mathbb{R}$ has been investigated leading to the time-independent problem 
\begin{equation}\label{6.3}
-\Delta u+(m_0^2-\omega^2)u-|u|^{p-2}u=0
\end{equation}
See also \cite{WalterA}. Already with the familiar cubic-cubic case, numerical solutions have been developed in \cite{Shenggao-Zhou.Xiao-liang.Cheng-Showmore} for the one-dimensional system
\begin{equation}\label{8.1}
\left\{
\begin{array}{l}
i\frac{\partial u}{\partial t}+\frac{1}{2}\frac{\partial^2u}{\partial x^2}+(|u|^2+\beta|v|^2)u=0,\\
i\frac{\partial v}{\partial t}+\frac{1}{2}\frac{\partial^2v}{\partial x^2}+(\beta|u|^2+|v|^2)v=0,\\
u(x,0)=u_0(x)~~~~~~v(x,0)=v_0(x).
\end{array}
\right.
\end{equation}
on $\mathbb{R}$, where $u$ and $v$ stands for the wave amplitude considered in two polarizations. The parameter $\beta$ is related to phase modulation. $u_0$ and $v_0$ are fixed functions assumed to be sufficiently regular. 

Coupled NLS system has been analyzed for symmetries and exact solutions in \cite{Liu-Ping.Sen-Yue-Lou}. The problem studied is related to atmospheric gravity waves governed by the following general coupled NLS system.
\begin{equation}\label{10.2}
\left\{
\begin{array}{l}
iu_t+\alpha_1u_{xx}+(\sigma_1|u|^2+\Gamma_1|v|^2)u=0,\\
iv_t+iCv_{x}+\alpha_2v_{xx}+(\Gamma_2|u|^2+\sigma_2|v|^2)v=0.
\end{array}
\right.
\end{equation}
It is noticed that such a problem may be transformed to the wall known Boussinesq equation. 

In \cite{H-AMINIKHAH.F-POURNASIRI.F-MEHRDOUST},  novel effective approach for systems of coupled NLS equations has been developed for the model problem
\begin{equation}\label{11.1}
\left\{
\begin{array}{lll}
iu_t+iu_x+u_{xx}+u+v+\sigma_1f(u,v)u=0,\\
iv_t-iv_x+v_{xx}+u-v+\sigma_2g(u,v)v=0,
\end{array}
\right.
\end{equation}
where $f$ and $g$ are smooth nonlinear real functions depending on $(|u|^2,|v|^2)$ and $\sigma_1$, $\sigma_2$ are parameters. In \cite{Zehra-Pinara.Ekin-Delikta}, the folowing nonlinear cubic-quintic and coupling quintic system of NLS equations has been examined,
\begin{equation}\label{9.1}
\left\{\begin{array}{lll}
u_t+m_1u_{xx}=(\alpha+i\beta)u+f_1|u|^2u+f_2|u|^4u+f_3|v|^2u+f_4|u|^2|v|^2u,\\
v_t+m_2v_{xx}=(\gamma+i\delta)v+g_1|v|^2v+g_2|v|^4v+g_3|u|^2v+g_4|u|^2|v|^2v,\\
\end{array}
\right.
\end{equation}
where $u$ and $v$ describe the complex envelopes of an electric field in a co-moving frame. $\alpha$ and $\gamma$ are potentials. $\beta$ and $\delta$ are amplifications. The functions $f_j$ and $g_j$, $j=1,1,3,4$ describe the cubic, quintic and coupling quintic nonlinearities coefficients respectively. $m_1$ and $m_2$ stands for the dispersion parameters

In \cite{HAIDONG}, existence of ground state solutions has been studied for the NLS system 
\begin{equation}\label{17.1}
\left\{
\begin{array}{lll}
-i\frac{\partial}{\partial t}\psi_1=\Delta\psi_1-v_1(x)\psi_1+\mu_1|\psi_1|^2\psi_1+\beta|\psi_2|^2\psi_1+\gamma\psi_2,\\
-i\frac{\partial}{\partial t}\psi_2=\Delta\psi_2-v_2(x)\psi_2+\mu_2|\psi_1|^2\psi_2+\beta|\psi_1|^2\psi_2+\gamma\psi_1,\\
\psi_j=\psi_j(x,t)\in\mathbb{C},
\end{array}
\right.
\end{equation}
$x\in\mathbb{R}$, $t>0$, $j=1,2$.

In \cite{T.Saanouni} a multi-nonlinearities coupled focusing NLS system has been studied for existence of ground state solutions and global existence and finite-time blow-up solutions. The authors considered precisely the coupled system 
\begin{equation}\label{24.1}
\left\{
\begin{array}{l}
i \dot{u}_j +\Delta u_j =-\displaystyle\sum^m_{k=1}a_{jk}|u_k|^{p-2}u_j,\\
u_j(0,.)=\psi_j
\end{array}
\right.
\end{equation}
where $u_j:\mathbb{R}^N \times \mathbb{R}\rightarrow\mathbb{C}$, $j=1,2,\dots,m$ and the $a_{jk}$'s are positive real numbers with $a_{jk}=a_{kj}$. 

In the present work, we focus on the nonlinear mixed super-linear cubic defocusing model  
\begin{equation}\label{f-modela}
f_1(u,v)=g(u,v)u=(|u|^{p-1}+\lambda|v|^{2})u=f_2(v,u),
\end{equation}
with $\lambda>0$ and $1<p\not=3$. We consider the evolutive nonlinear Schr\"odinger system on $\mathbb{R}\times(0,+\infty)$,
\begin{equation}\label{ContinuousProblemNLS1}
\left\{\begin{array}{lll}
iu_{t}+\sigma_1u_{xx}+g(u,v)u=0,\\
iv_{t}+\sigma_2v_{xx}+g(v,u)v=0.
\end{array}\right.
\end{equation}
Focuses will be on the study of the behavior of the steady state solutions according to their initial values. We propose precisely to develop a classification of the steady state solutions of problem (\ref{ContinuousProblemNLS1}). Recall that a steady state solution of problem (\ref{ContinuousProblemNLS1}) is any solution of the form
$$
W(x,t)=(e^{i\omega t}u(x),e^{i\omega t}v(x)),
$$
$\omega>0$. Substituting $W$ in problem (\ref{ContinuousProblemNLS1}), we immediately obtain a solution $(u,v)$ of the system
\begin{equation}\label{ContinuousSteadyStateProblem1-Dim1}
\begin{cases}
\sigma_1u_{xx}-\omega u+g(u,v)u=0,\\
\sigma_2v_{xx}-\omega v+g(v,u)v=0.
\end{cases}
\end{equation}
We will see that classifying the solutions of problem (\ref{ContinuousSteadyStateProblem1-Dim1}) depends strongly on the positoive/negative/null zones of the nonlinear function  
$$
g_\omega(x,y)=|x|^{p-1}+\lambda y^2-\omega,\;\;(x,y)\in\mathbb{R}^2
$$
which in turns varies accordingly to the parameters $p$, $\lambda$ and $\omega$. Denote also for $s\in\mathbb{R}$, $s\not=1$ and $\eta>0$, $\omega_{s,\eta}=\left(\displaystyle\frac{\omega}{\eta}\right)^{\frac{1}{s-1}}$. Denote also
$$
\Gamma^{\omega}_{p,\lambda}=\left\{(u,v)\in\mathbb{R}^2\,;\;|u|^{p-1}+\lambda v^2-\omega=0\right\}
$$
and
$$
\Gamma^{\omega,*}_{p,\lambda}=\left\{(u,v)\in\mathbb{R}^2\,;\;|v|^{p-1}+\lambda u^2-\omega=0\right\}.
$$
It is noticeable that such curves are more and more smooth whenever the parameter $p$ increases. To illustrate this fact, we provided in Appendix \ref{2emeAppendix} a brief overvew for some cases of $\Gamma^{\omega}_{p,\lambda}$ and $\Gamma^{\omega,*}_{p,\lambda}$. 

The paper is organized as follows. The next section is concerned with the development of our main results. Concluding and future directions are next raised in section 3. An appendix is also developed to illustrate the effects of the parameters of the problem on the classification of the solutions.
\section{Main results}
\begin{figure}[htbp]
	\begin{center}
		\includegraphics[scale=0.85]{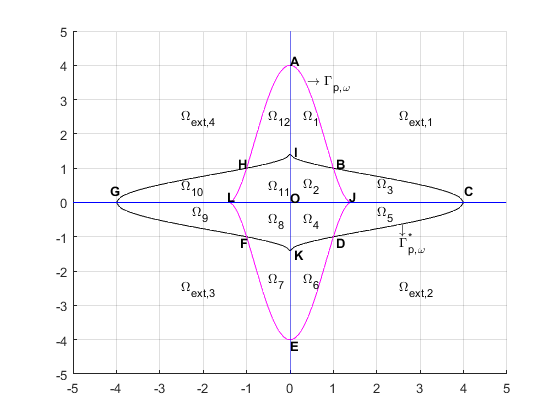}
		\caption{Partition of the plane $\mathbb{R}^2$ according to the curves $\Gamma_{p,\omega}$ and $\Gamma^{*}_{p,\omega}$ for $p=1.25$ and $\omega=2$.
		}\label{PartitionFigure1}
	\end{center}
\end{figure}

\begin{figure}[htbp]
	\begin{center}
		\includegraphics[scale=0.65]{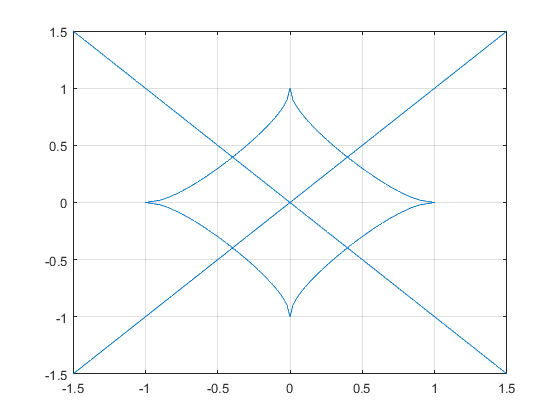}
		\caption{Partition of the plane $\mathbb{R}^2$ according to the curve $\Lambda_{p}$ (in blue) for $p=1.5$.}\label{PartitionFigure2}
	\end{center}
\end{figure}

\begin{figure}[htbp]
	\begin{center}
		\includegraphics[scale=0.65]{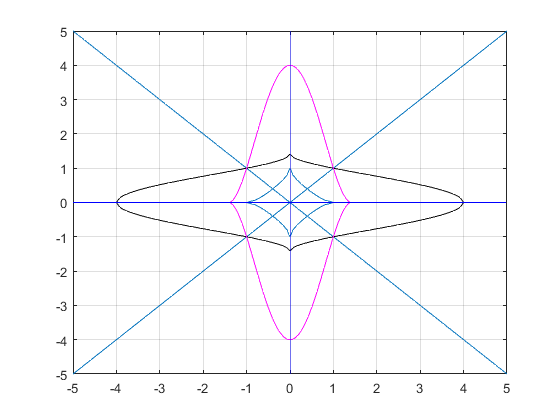}
		\caption{Partition of the plane $\mathbb{R}^2$ according to $\Lambda_{p,\omega}$, $\Gamma^{*}_{p,\omega}$ and $\Lambda_{p}$ for $p=1.5$ and $\omega=2$.}\label{PartitionFigure3}
	\end{center}
\end{figure}

As it is said in the introduction, the behavior of the solutions depends strongly on the parameters of the problem, 
especially those affacting the sign of the function $g_{\omega}$. It holds in fact that some of these parameters may be reduced to 
$$
\sigma_1=\sigma_2=\lambda=1,
$$
which simplifies the computations needed later. Indeed, denote
$$
u(r)=\widetilde{u}(\displaystyle\frac{r}{\sqrt\sigma_1})\quad\hbox{and}\quad v(r)=\widetilde{v}(\displaystyle\frac{r}{\sqrt\sigma_2}),
$$
the functions $\widetilde{u}$ and $\widetilde{v}$ satisfy immediately the system
$$
{u}_{xx}-\omega{u}+g({u},{v}){u}=0,\qquad
{v}_{xx}-\omega{v}+g({v},{u}){v}=0.
$$
Moreover, consider the scaling modifications
$$
u(r)=K_1\overline{u}(\alpha r)\quad\hbox{and}\quad v(r)=K_2\overline{v}(\beta r),
$$
where $K_1$, $K_2$, $\alpha$ and $\beta$ are constants to be fixed. The functions $\overline{u}$ and $\overline{v}$ satisfy the system
\begin{equation}\label{ContinuousProblemNLS1-a}
\begin{cases}
u_{xx}+(|u|^{p-1}+|v|^{2}-\omega)u=0,\\
v_{xx}+(|v|^{p-1}+|u|^{2}-\omega)v=0,\\
\end{cases}
\end{equation}
whenever the constants $K_1$, $K_2$, $\alpha$ and $\beta$ satisfy
\begin{equation}\label{ContinuousProblemNLS1-ConstantsSystem}
K_1^{p-1}=\sigma_1\alpha^2,\;\;K_2^{p-1}=\sigma_2\beta^2,\;\;
\lambda K_1^{2}=\sigma_2\beta^2,\;\;\lambda K_2^{2}=\sigma_1\alpha^2,
\end{equation}
which in turns yields that
$$
\alpha=\exp\bigl(\displaystyle\frac{(p-1)B_\lambda(\sigma_1,\sigma_2)-2A_\lambda(\sigma_1,\sigma_2)}{(p-3)(p+1)}\bigr),
$$
$$
\beta=\exp\bigl(\displaystyle\frac{2B_\lambda(\sigma_1,\sigma_2)-(p-1)A_\lambda(\sigma_1,\sigma_2)}{(p-3)(p+1)}\bigr),
$$
$$
K_1=\sqrt{\displaystyle\frac{\sigma_2\beta^2}{\lambda}}\;\;\hbox{and}\;\;K_2=\sqrt{\displaystyle\frac{\sigma_1\alpha^2}{\lambda}},
$$
where
$$
A_\lambda(\sigma_1,\sigma_2)=\displaystyle\frac{1}{1-p}\log\sigma_1+\displaystyle\frac{1}{2}\log\sigma_2-\displaystyle\frac{1}{2}\log\lambda
$$
and
$$
B_\lambda(\sigma_1,\sigma_2)=-\displaystyle\frac{1}{2}\log\sigma_1+\displaystyle\frac{1}{p-1}\log\sigma_2-\displaystyle\frac{1}{2}\log\lambda.
$$
Given these facts, we will consider in the rest of the paper the problem (\ref{ContinuousProblemNLS1-a}) with the initial conditions 
\begin{equation}\label{ContinuousProblemNLS1-b}
u(0)=a\,,\;v(0)=b\,,\;\;u'(0)=v'(0)=0,
\end{equation}
where $a,b\in\mathbb{R}$. The associated curves will be denoted simply 
$$
\Gamma_{p,\omega}=\left\{(u,v)\in\mathbb{R}^2\,;\;|u|^{p-1}+|v|^2-\omega=0\right\}
$$
and
$$
\Gamma^{*}_{p,\omega}=\left\{(u,v)\in\mathbb{R}^2\,;\;|v|^{p-1}+|u|^2-\omega=0\right\}
$$
and
$$
\Lambda_{p}=\left\{(u,v)\in\mathbb{R}^2\,;\;|u|^{p-1}+|v|^2=|v|^{p-1}+|u|^2\right\}.
$$
Denote also $\omega_{s}=\omega^{\frac{1}{s-1}}$, for $s\in\mathbb{R}$, $s\not=1$. Figure \ref{PartitionFigure1} illustrates the partition of the plane $\mathbb{R}^2$ according to the curves $\Gamma_{p,\omega}$ and $\Gamma^{*}_{p,\omega}$. Figure \ref{PartitionFigure2} illustrates the partition of the plane $\mathbb{R}^2$ according to the curve $\Lambda_{p}$. Finally, Figure \ref{PartitionFigure3} illustrates the partition of the plane $\mathbb{R}^2$ according to all the curves $\Lambda_{p,\omega}$, $\Gamma^{*}_{p,\omega}$ and $\Lambda_{p}$.
\\
\\

We now start developing our main results. For this aim, we assume in the rest of the paper that
\begin{equation}\label{assumption-on-w-and-p}
\omega\geq1\;\;\hbox{and}\;\;1<p<3.
\end{equation}
It is straightforward that $\omega_{p}>\omega_{3}$. Moreover, the closed curves $\Gamma_{p,\omega}$ and $\Gamma^*_{p,\omega}$ intersect at four points $B$, $D$, $F$ and $H$ as it is shown in Figure \ref{PartitionFigure1}. Notice also that the curves $\Gamma_{p,\omega}$ and $\Gamma^*_{p,\omega}$ split the plane $(u,v)$ into twelve regions, $\Omega_i$, $i=1,2,...,12$ and an external region $\Omega_{ext}=\displaystyle\bigcup_{1\leq i\leq4}\Omega_{ext,i}$. Such regions are in the heart of the classification of solutions of problem
(\ref{ContinuousProblemNLS1-a})-(\ref{ContinuousProblemNLS1-b}). Notice also from the symmetry of the function $g_\omega$ that the points $B$ and $F$ satisfy the cartezian equation $u=v$ and the points $D$ and $H$ satisfy $u=-v$. Furthermore the polygon $BDFH$ is a square.

This symmetry also shows easily that whenever $(u,v)$ is a pair of solutions of the system (\ref{ContinuousProblemNLS1-a}), the pairs $(-u,v)$, $(u,-v)$ and $(-u,-v)$ are also solutions. Consequently, in the rest of the paper we will focus only on the case where $u(0)=a\geq0$ and $v(0)=b\geq0$. The remaining cases will be deduced by symmetry. As the function $f(u,v)$ is Locally Lipshitz continues, the existence and uniqueness of the solution is guaranteed by means of the famous Cauchy-Lipschitz theory. This  will be recalled in the Appendix later.
\begin{theorem}\label{abinOmega1}
	Whenever the initial data $(u(0),v(0))=(a,b)\in\Omega_1$, the problem (\ref{ContinuousProblemNLS1-a})-(\ref{ContinuousProblemNLS1-b}) has a unique positive solution lieing in the same region $\Omega_1$. Furthermore,
	\begin{description}
		\item[i.] $u$ is nondecreasing,
		\item[ii.] $v$ is nonincreasing,
		\item[iii.] $(u,v)$ tends to $B$ as $x$ tends to $\infty$.
	\end{description}
\end{theorem}
\textbf{Proof.}
	Writing problem (\ref{ContinuousSteadyStateProblem1-Dim1}) at $x=0$, yields that
	$$
	u''(0)=-g_\omega(a,b)a>0\qquad\mbox{and}\qquad v''(0)=-g_\omega(b,a)b<0.
	$$
	Hence, there exists $\delta>0$ small enought such that
	$$
	u''(x)>0\qquad\mbox{and}\qquad v''(x)<0,\;\;\forall x\in(-\delta,\delta).
	$$
	Consequently, $u'$ is nondecreasing on $(-\delta,\delta)$ and $v'$ is nonincreasing on $(-\delta,\delta)$. Next, as $u'(0)=v'(0)=0$, it therefore follows that $u$ is nonincreasing on $(-\delta,0)$ and nondecreasing on $(0,\delta)$ and that $v$ is nondecreasing on $(-\delta,0)$ and nonincreasing on $(0,\delta)$. Assume that $u$ remains nondeceasing on $(0,\infty)$. Then, we get at a first step $u(x)>a$, $\forall x>0$. Two situations may occur.
	$u$ has a finite limit $l_u$ which is obviously positive, or $u$ is going up to $\infty$ as $x$ goes up to $\infty$. Whenever the last situation occurs, we immediately get
	$$
	u(x)\leq-x^2+K_1x+K_2
	$$
	for some constants $K_1$ and $K_2$ and for $x$ large enough. Consequently, $u(x)\rightarrow-\infty$ as $x\rightarrow+\infty$. Which is a contradiction. Now, if the first case occurs, then, we immediately get a limit $L_v\geq0$ of $v^2(x)$ as $x$ goes up to $\infty$. Thus, denoting $l_v=\sqrt{L_v}$, we get
	$$
	g_\omega(l_u,l_v) =g_\omega(l_v,l_u)=0.
	$$
	We claim that $(l_u,l_v)=B$. Indeed, if not, this means that the solution $(u,v)$ crosses one of the lines $(A,B)$ or $(B,I)$ at one point $X(x_0)$. At $x_0$, we have
	$$
	x_0>0,\;\;g_\omega(u(x_0),v(x_0))=0,\;\; g_\omega(v(x_0),u(x_0))>0.
	$$
	Consider next the energy functional $E$ associated to problem (\ref{ContinuousSteadyStateProblem1-Dim1}) and defined by
	$$
	E(u,v)(x)=\displaystyle\frac{1}{2}(u_x^2+v_x^2)+\displaystyle\frac{1}{p+1}(|u|^{p+1}+|v|^{p+1})+\displaystyle\frac{\lambda}{2}u^2v^2-\displaystyle\frac{\omega}{2}(u^2+v^2).
	$$
	It is straightforward that $E$ is constant as a function of $x$. Henceforth,
	$$
	E(u,v)(x)=E(u,v)(0)=E(u,v)(\infty).
	$$
	This equation may be written otherwise as
	$$
	\displaystyle\frac{1}{p+1}(a^{p+1}+b^{p+1})+\displaystyle\frac{\lambda}{2}a^2b^2-\displaystyle\frac{\omega}{2}(a^2+b^2)=
	\displaystyle\frac{2}{p+1}l^{p+1}+\displaystyle\frac{\lambda}{2}l^4-\omega l^2,
	$$
	where $l=l_u=|l_v|$ is the common limit of $u$ and $|v|$ at $\infty$. Denote next
	$$
	H(u,v)=\displaystyle\frac{1}{p+1}(|u|^{p+1}+|v|^{p+1})+\displaystyle\frac{\lambda}{2}u^2v^2-\displaystyle\frac{\omega}{2}(u^2+v^2).
	$$
	Whenever $u\in(a,u(x_0))$ and $v\in(v(x_0),b)$ or otherwise $x\in(0,x_0)$, the functional $H$ is non-increasing as a function of $u$ and non-decreasing as a function of $v$. Hence forth,
	$$
	E(u,v)(x_0)<E(u,v)(x)<E(u,v)(0),
	$$
	which is contradictory. As a consequence, we can state that $(u,v)$ does not cross the arc $(A,B)$. In addition, by similar techniques we can show that $(u,v)$ does not cross the arc $(I,B)$. So, necessarily, $(u,v)$ remains bounded in the domain $\Omega_1$. This yields that, $(u,v)$ tends to $B$ as $x$ tends to $\infty$.

\begin{theorem}\label{abinOmega2}
	Whenever the initial data $(u(0),v(0))=(a,b)\in\Omega_2$, the problem (\ref{ContinuousProblemNLS1-a})-(\ref{ContinuousProblemNLS1-b}) has a unique positive solution lieing in the same region $\Omega_2$. Furthermore,
	\begin{description}
		\item[i.] $u$ and $v$ are nondecreasing,
		\item[ii.] $(u,v)$ tends to $B$ as $x$ tends to $\infty$.
	\end{description}
\end{theorem}
\textbf{Proof.}
	We proceed as in the proof of Theorem \ref{abinOmega1}. At $x=0$, we have
	$$
	u''(0)=-g_\omega(a,b)a>0\qquad\mbox{and}\qquad v''(0)=-g_\omega(b,a)b>0.
	$$
	Therefore, as $u'(0)=v'(0)=0$, it follows that $u$ and $v$ are nondecreasing on $(0,\delta)$ for some $\delta>0$ small enough. Assume that $u$ remains nondecreasing on $(0,\infty)$. Then, we get at a first step $u(x)>a$, $\forall x>0$. Two situations may occur. $u$ has a finite limit $l_u$ which is obviously positive, or $u$ is going up to $\infty$ as $x$ goes up to $\infty$. Proceeding as in the proof of Theorem \ref{abinOmega1}, we show that the second situation can not occur here also. So, assume that the second the first case occurs. Then, we immediately get a limit $L_v\geq0$ of $v^2(x)$ as $x$ goes up to $\infty$. Denote as previously $l_v=\sqrt{L_v}$. Whenever $l_v=0$ (or equivalently $L_v=0$), then $l_u=\omega_p$. Consequently, $(u,v)\longrightarrow\,J$ as $x\longrightarrow+\infty$. Therefore the energy $E$ at $\infty$ becomes
	$$
	E(l_u,0)=l_u^2\omega\displaystyle\frac{1-p}{2(1+p)}>0.
	$$
	On the other hand
	$$
	E(u,v)(0)=\displaystyle\frac{a^{p+1}+b^{p+1}}{p+1}+\displaystyle\frac{\lambda}{20a^2b^2}-\displaystyle\frac{\omega}{2}(a^2+b^2)>0.
	$$
	So, we get a contradiction. Hence, $v$ is bounded but does not tend to 0 at $\infty$. Consequently, as $v^2$ has a limit $L_v$ at $\infty$, its derivative tends to 0 at $\infty$. So, necessarily $v'$ and consequently $v''$ tend to 0 at $\infty$. The second equation of problem (\ref{ContinuousSteadyStateProblem1-Dim1}) yields then that $g_\omega(l_v,l_u)=0$. On $\partial\Omega_2$ it remains only the point $B$ to be the candidate for $(l_v,l_u)$. \\
	Assume now that $u$ is not monotone on $(0,\infty)$ and let $x_0$ be its the first positive critical point. We get
	$$
	E(u,v)(x_0)>\displaystyle\frac{u(x_0)^{p+1}+v(x_0)^{p+1}}{p+1}+\displaystyle\frac{\lambda}{2}u(x_0)^2v(x_0)^2-\displaystyle\frac{\omega}{2}(u(x_0)^2+v(x_0)^2)>0.
	$$
	Whenever $(u(x_0),v(x_0))\in\Omega_2$, this leads to
	$$
	E(u,v)(x_0)>E(u,v)(0).
	$$
	Which is contradictory. Next, whenever $(u(x_0),v(x_0))\notin\Omega_2$, consider the point $x_1$ to be the point at which $(u,v)$ crosses $\partial\Omega_2$. Three positions may hold. $(u(x_1),v(x_1))\in(B,J)$ or $(u(x_1),v(x_1))\in(I,B)$ or $(u(x_1),v(x_1))=B$. The last one leads to a stationary solution $(u,v)=B$. In the first case, we get again
	$$
	E(u,v)(x_1)>E(u,v)(0),
	$$
	which is a contradiction (Similarly, for the second case). So, to conclude, we have proved that $(u,v)$ does not cross $\overline{\Omega_2}$. Thus the unique attractive point is $B$.

Now, similarly to the previous cases, we may prove the following result.
\begin{theorem}\label{abinOmega3}
	Whenever the initial data $(u(0),v(0))=(a,b)\in\Omega_3$, the problem (\ref{ContinuousProblemNLS1-a})-(\ref{ContinuousProblemNLS1-b}) has a unique positive solution lieing in the same region $\Omega_3$. Furthermore,
	\begin{description}
		\item[i.] $u$ is nonincreasing,
		\item[ii.] $v$ is nondecreasing,
		\item[iii.] $(u,v)$ tends to $B$ as $x$ tends to $\infty$.
	\end{description}
\end{theorem}
We now study the case when the initial value is out of the negative zones of $g_\omega(u,v)$ and $g_\omega(v,u)$. We have the following result.
\begin{theorem}\label{abinOmegaext}
	Whenever the initial data $(u(0),v(0))=(a,b)\in\Omega_{ext}$, the problem (\ref{ContinuousProblemNLS1-a})-(\ref{ContinuousProblemNLS1-b}) has a unique solution. Furthermore,
	$u$ and $v$ are not simultaneously monotone and satisfy
	$$
	\lim_{x\rightarrow\infty}(u,v)(x)\in\Gamma^{\omega}_{p,\lambda}\cap\Gamma^{\omega,*}_{p,\lambda}.
	$$
\end{theorem}
\textbf{Proof.}
	We will prove the theorem for the case $(a,b)\in\Omega_{ext,1}$. The remaining ones may be proved by applying similar techniques. The starting point is always the same. We evaluate the behaviour of the solution at the origin $x=0$. We have
	$$
	u''(0)=-g_\omega(a,b)a<0\qquad\mbox{and}\qquad v''(0)=-g_\omega(b,a)b<0.
	$$
	Hence, as $u'(0)=v'(0)=0$, $u$ and $v$ are nonincreasing on $(0,\delta)$ for some $\delta>0$ small enough. Assume that $u$ remains nonincreasing on $(0,\infty)$. Two situations may occur. $u$ has a finite limit $l_u$, or $u$ is going down to $-\infty$ as $x$ goes up to $\infty$. Proceeding as in the proof of Theorem \ref{abinOmega1}, we show that the second situation can not occur here-also. So, assume that the first case occurs. Then, we immediately get a limit $L_v\geq0$ of $v^2(x)$ as $x$ goes up to $\infty$. Denote as previously $l_v=\sqrt{L_v}$. It is straightforward that $g_\omega(l_u,l_v)l_u=0$. We claim next that $g_\omega(l_v,l_u)=0$, Indeed, if it is not, we get 
	$$
	v''(x)\simeq-g_\omega(l_v,l_u)v\,;\,\,x\rightarrow\infty.
	$$
	This means that 
	$$
	v(x)=K_1\cos(\omega x+\varphi)\,;\,\,x\rightarrow\infty
	$$
	or 
	$$
	v(x)=K_1e^{-\omega x}+K_2e^{\omega x}\,;\,\,x\rightarrow\infty,
	$$
	for some real constants $K_1$, $K_2$ and $\omega>0$. In the first case, $v^2$ does not have a limit as $x\rightarrow\infty$, except that $K_1=0$ and thus $v\equiv0$ for $x$ large enough. In the second case, we get firstly $K_2=0$ for the same reason on $v^2$. Therefore, $v(x)\rightarrow0$ as $x\rightarrow\infty$. So that, $l_u=0$ also. Which is contradictory. So now, we have obviously, $l_u=0$ or $l_u\neq0$. Whenever the latter occurs we get $g_\omega(l_u,l_v)=g_\omega(l_v,l_u)=0$ and $l_ul_v\neq0$. Which means that $(l_u,l_v)=\Gamma^\omega_{p,\lambda}\cap\Gamma_{p,\lambda}^{\omega,*}$. So we get furthemore $l_u=l_v=l$, where  
	$$
	l^{p-1} +\lambda l^2-w=0.
	$$
	Now, applying the energy functional $E(u,v)$ and $H(u,v)$ we get
	$$
	\displaystyle\frac{\partial H}{\partial u}=|u|^p+\lambda u v^2-wu=u(|u|^{p-1})+\lambda v^2-w=uf(u,v)>0
	$$
	and
	$$
	\displaystyle\frac{\partial H}{\partial v}=|v|^p+\lambda u^2 v-wu=u(|u|^{p-1})+\lambda v^2-w= vf(v,u)>0.
	$$
	Consequently,
	$$
	H(u,v)(0)=E(u,v)(0)< H(u,v)(\infty)=E(u,v) (\infty)
	$$
	which is contradictory. Consequently, one of the  following situations hold.
	\begin{itemize}
		\item $u\searrow$, $v\searrow$, $(u,v)\longrightarrow B$ and $u\geq x_B$, $v\geq y_B$.
		\item $u$ and $v$ are not monotone on the whole interval $(0,\infty)$ and thus $(u,v)$ crosses the boundary $\partial\Omega^c_{ext}$ to $\Omega_1$ or $\Omega_2$ and thus $(u,v)$ is attracted by $B$ again. 
	\end{itemize}
	Now, we will study the case $\ell_u=0$. We get in this case
	$$
	v''+(|v|^{p-1}+\lambda u^2)v=0.
	$$
	At $\infty$ we obtain 
	$$
	v''+|v|^{p-1}v=0.
	$$
	$v$ is $\searrow$ on $(0,\infty)$ for $\delta>0$ small enough. If it remains $\searrow$ on $(0,\infty$), it has a limit $l_v$ as $x\longrightarrow+\infty$. Whenever $l_v=-\infty$ we get for $x\longrightarrow +\infty$, 
	$$
	v(x)\geq x^2+k_1x+k_2,
	$$
	for some constants $k_1,k_2$, which is contradictory. Whenever $l_v$ is finite and nonzero, we get at infinity
	$$
	v(x)=k \cos (\delta x+\varphi)
	$$
	which has no limit as $x\longrightarrow +\infty$. Consequently, $v$ can not be $\searrow$ on $(0,\infty)$. It remains that $u\searrow 0$ and $v\searrow 0$. This way, $(u,v)$ crosses $\Omega_{ext}$ to $\Omega_{int}$. Thus $(u,v)$ will be attracted by $B$ which contradicts the monotonicity.
	
	We conclude thus that $u$ and $v$ can not be monotone simultaneously and that $(u,v)$ tends to one of the points $B,D,F,H$ of $\Gamma^\omega_{p,\lambda}\cap \Gamma_{p,\lambda}^{\omega,*}$.
	
We now study the case when the initial value lies on one of the curves $\Gamma_{p,\lambda}^\omega\cup\Gamma^{\omega,*}_{p,\lambda}$. Of course, whenever the initial value $(u(0),v(0))\in\Gamma_{p,\lambda}^\omega\cap\Gamma^{\omega,*}_{p,\lambda}$, the pair solution $(u,v)$ remains constant due to Cauchy-Lipschitz theorem. The interesting cases will be those when the initial value is on the graph $\Gamma_{p,\lambda}^\omega\cup\Gamma^{\omega,*}_{p,\lambda}\setminus\{B,D,F,H\}$. It holds in fact that the remaining vertices of the graph necessitates to be treated case by case. For this aim and throughout the remaining parts of the present paper, we denote $G$ the graph $\Gamma_{p,\lambda}^\omega\cup\Gamma^{\omega,*}_{p,\lambda}$ and $G^\star$ the graph $G$ without vertices. The following result is obtained.
\begin{theorem}\label{absurlegrapheGstar}
	Whenever the initial data $(u(0),v(0))=(a,b)\in G^\star$, the problem (\ref{ContinuousProblemNLS1-a})-(\ref{ContinuousProblemNLS1-b}) has a unique solution. Furthermore,
	$u$ and $v$ are not simultaneously monotone and satisfy
	$$
	\lim_{x\rightarrow\infty}(u,v)(x)\in\Gamma_{p,\lambda}^\omega\cap\Gamma_{p,\lambda}^{\omega,\star}.
	$$
\end{theorem}
\textbf{Proof.}
	At the origin $x=0$, we get 
	$$
	u''(0)=u'(0)=-g_\omega(a,b)a=0\qquad\mbox{and}\qquad v''(0)=-g_\omega(b,a)b<0.
	$$
	We claim that $u$ can not remain constant on any interval $(0,\delta)$, $\delta>0$. Indeed, if this occurs on some interval $(0,\delta)$, we get immediately 
	$$
	u(x)=a\quad\mbox{and}\quad\lambda v^2(x)=\omega-a^{p-1}\,,\;\;\forall\,x\in(0,\delta).
	$$
	So, $|v|$ is also constant ($\equiv b$) on $(0,\delta)$. Therefore, we get
	$$
	v''+g_\omega(b,a)v=0\;\;\mbox{on}\;\;(0,\delta).
	$$
	This means that
	$$
	v(x)=K_1\cos(\rho x+\rho)\,;\,\,x\in(0,\delta),
	$$
	which contradicts the fact that $|v|$ is constant, except if $K_1=0$. However, in this case we again get $b=0$, which is contradictory. So as the claim. Consequently, three situations may occur. $(i)$ the pair solution $(u,v)$ remains on the edge $)A,B($ and being non constant. $(ii)$ the pair solution $(u,v)$ bifurcates upward $\Omega_{ext,1}$. $(iii)$ the pair solution $(u,v)$ bifurcates toward $\Omega_{1}$.
	The case (i): $(u,v)$ remains on $)A,B($. We get $u(x)=a, \forall x \in (0, \infty)$ which is a contradictory with the previous claim. Assume now that the situation (ii) or (iii) occur. Thus it is attracted by the limit point $B$. 

Next, observing that the nonlinear function $f(u,v)=(|u|^{p-1}+\lambda|v|^2)u$ has some parity characteristics relatively to $u$ and $v$, the analogous of Theorems \ref{abinOmega1}, \ref{abinOmega2}, \ref{abinOmega3}, \ref{abinOmegaext}, \ref{absurlegrapheGstar} may be obtained. We precisely get
\begin{theorem}\label{abinOmega4-12}
	Whenever the initial data $(u(0),v(0))=(a,b)\in\displaystyle\bigcup_{i=4}^{12}\Omega_i$, the problem (\ref{ContinuousProblemNLS1-a})-(\ref{ContinuousProblemNLS1-b}) has a unique solution lieing in the same region $\Omega_i$ that contains its initial value $(a,b)$ . Furthermore,
	\begin{enumerate}
		\item For $(a,b)\in\Omega_4$,
		\begin{description}
			\item[i.] $u$ is positive and nonincreasing,
			\item[ii.] $v$ is negative and nondecreasing,
			\item[iii.] $(u,v)$ tends to $D$ as $x$ tends to $\infty$.
		\end{description}
		\item For $(a,b)\in\Omega_5$,
		\begin{description}
			\item[i.] $u$ is positive and nondecreasing,
			\item[ii.] $v$ is negative and nondecreasing,
			\item[iii.] $(u,v)$ tends to $D$ as $x$ tends to $\infty$.
		\end{description}
		\item For $(a,b)\in\Omega_6$,
		\begin{description}
			\item[i.] $u$ is positive and nondecreasing,
			\item[ii.] $v$ is negative and nondecreasing,
			\item[iii.] $(u,v)$ tends to $D$ as $x$ tends to $\infty$.
		\end{description}
		\item For $(a,b)\in\Omega_7$,
		\begin{description}
			\item[i.] $u$ is negative and nonincreasing,
			\item[ii.] $v$ is negative and nondecreasing,
			\item[iii.] $(u,v)$ tends to $F$ as $x$ tends to $\infty$.
		\end{description}
		\item For $(a,b)\in\Omega_8$,
		\begin{description}
			\item[i.] $u$ is negative and nonincreasing,
			\item[ii.] $v$ is negative and nonincreasing,
			\item[iii.] $(u,v)$ tends to $F$ as $x$ tends to $\infty$.
		\end{description}
		\item For $(a,b)\in\Omega_9$,
		\begin{description}
			\item[i.] $u$ is negative and nondecreasing,
			\item[ii.] $v$ is negative and nonincreasing,
			\item[iii.] $(u,v)$ tends to $F$ as $x$ tends to $\infty$.
		\end{description}
		\item For $(a,b)\in\Omega_{10}$,
		\begin{description}
			\item[i.] $u$ is negative and nondecreasing,
			\item[ii.] $v$ is positive and nondecreasing,
			\item[iii.] $(u,v)$ tends to $H$ as $x$ tends to $\infty$.
		\end{description}
		\item For $(a,b)\in\Omega_{11}$,
		\begin{description}
			\item[i.] $u$ is negative and nonincreasing,
			\item[ii.] $v$ is positive and nondecreasing,
			\item[iii.] $(u,v)$ tends to $H$ as $x$ tends to $\infty$.
		\end{description}
		\item For $(a,b)\in\Omega_{12}$,
		\begin{description}
			\item[i.] $u$ is negative and nonincreasing,
			\item[ii.] $v$ is positive and nonincreasing,
			\item[iii.] $(u,v)$ tends to $H$ as $x$ tends to $\infty$.
		\end{description}
	\end{enumerate}
\end{theorem}
\section{Conclusion}
In \cite{Yaping}, a system of coupled NLS and Heat equations has been considered for exponential stability in a torus region. The basic idea consists of transforming the system into one-dimensional coupled one by applying polar coordinates.

Inspired from the present work, in forthcomming works, we intend to consider analogoue study with 
\begin{itemize}
	\item The Heat operator 
	$$
	\mathcal{L}_i(u(x,t))=\mathcal{H}_i(u(x,t))=u_t-\sigma_i\Delta u
	$$ 
	leading to a nonlinear Heat system.
	\item The mixed Schr\"odinger-Heat operator 
	$$
	\mathcal{L}(u(,v))=(\mathcal{L}_1(u),\mathcal{L}_2(v))=(\mathcal{S}(u),\mathcal{H}(v))
	$$ 
	leading to a nonlinear coupled system of Schr\"odinger-Heat type. 
\end{itemize}

\section{Appendix}
\subsection{Appendix A}
Recall that in the previous sections, we applied for many times the well-known Cauchy Lipschitz theorem on the existence and uniqueness of solutions. In this section and for convenience, we will show that the generator function used is already locally Lipschitz continuous. 

Denote $\varphi=u'$ and $\psi=v'$. The system (\ref{ContinuousProblemNLS1-a})-(\ref{ContinuousProblemNLS1-b}) becomes
\begin{equation}
\left\{
\begin{array}{lll}
u' =\varphi\\
\varphi'=-(|u|^{p-1}+|v|^2)u,\\
v'=\psi\\
\psi'=-(|v|^{p-1}+|u|^2)v,\\
u(0)=a,\;v(0)=b,\;\varphi(0)=0,\;\psi(0)=0.
\end{array}
\right.
\end{equation}
Denoting $X$ the vector $X=(u,v,\varphi,\psi)^T$, where $^T$ stands for the transpose, we get
$$
X'=F(X),\;X(0)=(a,b,0,0)^T,
$$
where $F$ is the function defined by 
$$
F(x_1,x_2,x_3,x_4)=(x_2,-g_\omega(x_1,x_3)x_1,x_4,-g_\omega(x_3,x_1)x_3),\;\;(x_1,x_2,x_3,x_4)\in\mathbb{R}^4.
$$
\begin{lemma}
	$F$ is locally Lipschitz continuous on $\mathbb{R}^4$.
\end{lemma}
\textbf{Proof.} Let $\delta>0$ small enough and $X^0=(x_1^0,x_2^0,x_3^0,x_4^0)\in\mathbb{R}^4$ be fixed. For all $X,Y$ in the ball $B(X^0,\delta)$ we have
	$$
	\begin{array}{lll}
	\|F(X)-F(Y)\|_2^2&=&(x_2-y_2)^2+(x_4-y_4)^2\\
	&&+(g(x_1,x_3)x_1-g(y_1,y_3)y_1)^2\\
	&&+(g(x_3,x_1)x_3-g(y_3,y_1)y_3)^2.
	\end{array}
	$$
	We shall now evaluate the quantity $(g(x_1,x_3)x_1-g(y_1,y_3)y_1)^2$. Similar techniques will lead to $(g(x_3,x_1)x_3-g(y_3,y_1)y_3)^2$. We have
	$$
	\begin{array}{lll}
	&&|g(x_1,x_3)x_1-g(y_1,y_3)y_1)|\\
	&=&||x_1|^{p-1}x_1-|y_1|^{p-1}y_1+\lambda(|x_3|^{2}x_1-|y_3|^{2}y_1)|\\
	&\leq&||x_1|^{p-1}x_1-|y_1|^{p-1}y_1|+\lambda||x_3|^{2}x_1-|y_3|^{2}y_1|\\
	&\leq&C_1(p,X_0,\delta)|x_1-y_1|+\lambda\left[||x_3|^{2}-|y_3|^{2}||x_1|+|x_1-y_1||y_3|^{2}\right]\\
	&\leq&C_1(p,X_0,\delta)|x_1-y_1|+\lambda C_2(X_0,\delta)\left[||x_3-y_3|+|x_1-y_1|\right]\\
	&\leq&C_1(p,X_0,\delta)|x_1-y_1|+\lambda C_2(X_0,\delta)\left[||x_3-y_3|+|x_1-y_1|\right]\\
	&\leq&C(p,X_0,\delta,\lambda)\left[||x_3-y_3|+|x_1-y_1|\right],
	\end{array}
	$$
	where $C_1(p,X_0,\delta)>0$ is a constant depending only on $p$, $X_0$ and $\delta$. $C_2(X_0,\delta)>0$ is a constant depending only on $X_0$ and $\delta$. $C(p,X_0,\delta,\lambda)>0$ is a constant depending only on $p$, $X_0$, $\delta$ and $\lambda$. Therefore,
	$$
	\begin{array}{lll}
	|g(x_1,x_3)x_1-g(y_1,y_3)y_1)|^2&\leq&C(p,X_0,\delta,\lambda)^2\left[||x_3-y_3|+|x_1-y_1|\right]^2\\
	&\leq&2C(p,X_0,\delta,\lambda)^2\left[||x_3-y_3|^2+|x_1-y_1|^2\right].
	\end{array}
	$$
	Similarly, 
	$$
	|g(x_3,x_1)x_3-g(y_3,y_1)y_3)|^2\leq2\widetilde{C}(p,X_0,\delta,\lambda)^2\left[||x_3-y_3|^2+|x_1-y_1|^2\right],
	$$
	with some constant $\widetilde{C}(p,X_0,\delta,\lambda)>0$ analogue to ${C}(p,X_0,\delta,\lambda)$. Combining all these inequalities, we obtain
	$$
	\|F(X)-F(Y)\|_2^2\leq C(p,X_0,\delta,\lambda)\|X-Y\|_2^2,\;\;\forall\,X,Y\in\,B(X_0,\delta).
	$$

\subsection{Appendix B}\label{2emeAppendix}
In this part, we investigate the dependence of the different regions and curves $\Omega_i$, $i=1,\dots,12$, $\Omega_{ext,i}$, $i=1,\dots,4$ and the different curves $\Gamma^{\omega}_{p,\lambda}$, $\Gamma^{\omega,*}_{p,\lambda}$ and $\Lambda_{p}$. So denote $\Lambda_{p}^*=\Lambda_{p}\setminus\{(u,v)\in\mathbb{R}^2;\;|u|\not=|v|\}$. Denote also $\Omega_{int}=\displaystyle\bigcup_{i=1}^{12}\Omega_i$ and $\mathcal{A}_p$ the interior area countered by the curve $\Lambda_p^*$. 
\begin{lemma}\label{inclusionds[-,1]}
	The erea $\mathcal{A}_p\subset[-1,1]^2$.
\end{lemma}
Indde, denote for $(u,v)\in\mathcal{A}_p$
$$
u=r\cos\theta,\;\;v=r\sin\theta.
$$
It is staightforward that $r$ is maximum whenever $(u,v)$ is on the frontier $\Lambda_{p}^*$, which may be then governed by the polar equation
$$
r^{p-3}\biggl(|\cos\theta|^{p-1}-|\sin\theta|^{p-1}\biggr)=\cos(2\theta).
$$
Otherwise,
$$
r^{p-3}=\displaystyle\frac{\cos(2\theta)}{|\cos\theta|^{p-1}-|\sin\theta|^{p-1}}.
$$
Straightforward calculus yield that whenever $0<p<3$, the radius $r$ is extremum for $\theta\in\{0,\pm\frac{\pi}{2},\pi\}$, which gives the vertices $(\pm1,0)$, $(0,\pm1)$. In fact $r(_theta)$ may be extended at $\theta=\{\pm\frac{\pi}{4},\pm\frac{3\pi}{4}\}$. However, these points yield immediately $|u|=|v|$, which is the trivial (unbounded) part of $\Lambda_{p}$.

As a result of Lemma \ref{inclusionds[-,1]}, we immediately deduce that whenever $0<\omega\leq1$ and $1<p<3$ the inclusion $\Omega_{int}\subset\mathcal{A}_p$. However, for $\omega\geq1$ and $1<p<3$, the area $\subset\mathcal{A}_p$ is contained in the curved octagonal shape $IBJDKFLH$.
\begin{figure}[htbp]
	\begin{center}
		\includegraphics[scale=0.65]{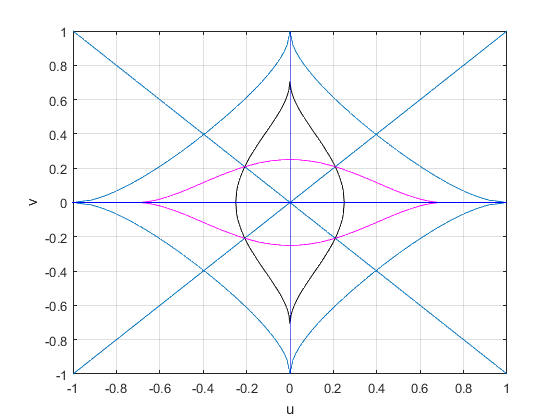}
		\caption{Parameters' domains for $p=1.5$, $\omega=0.5$}\label{w05p15}
	\end{center}
\end{figure}
\\
\begin{figure}[htbp]
	\begin{center}
		\includegraphics[scale=0.65]{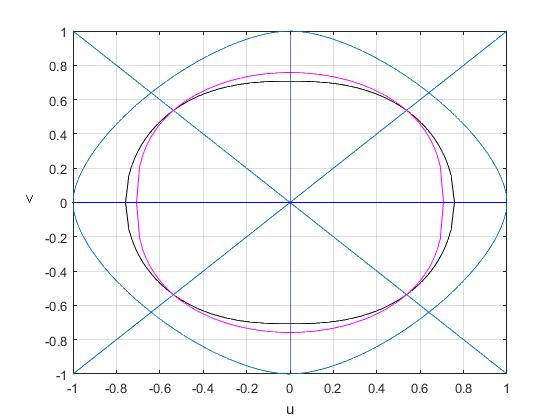}
		\caption{Parameters' domains for $p=3.5$, $\omega=0.5$}\label{w05p35}
	\end{center}
\end{figure}
\\
\begin{figure}[htbp]
	\begin{center}
		\includegraphics[scale=0.65]{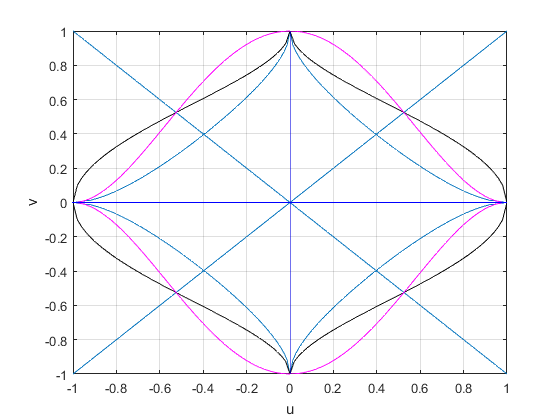}
		\caption{Parameters' domains for $p=1.5$, $\omega=1$}\label{w1p15}
	\end{center}
\end{figure}
\\
\begin{figure}[htbp]
	\begin{center}
		\includegraphics[scale=0.65]{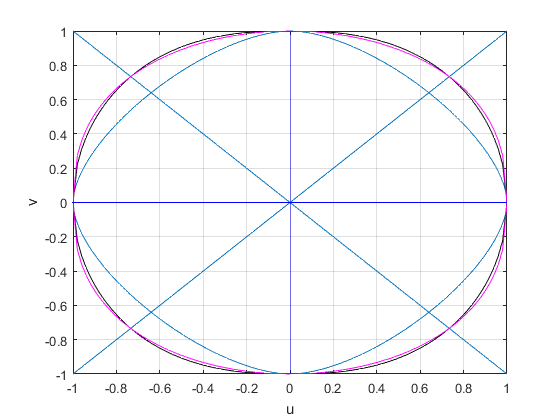}
		\caption{Parameters' domains for $p=3.5$, $\omega=1$}\label{w1p35}
	\end{center}
\end{figure}
\\
\begin{figure}[htbp]
	\begin{center}
		\includegraphics[scale=0.65]{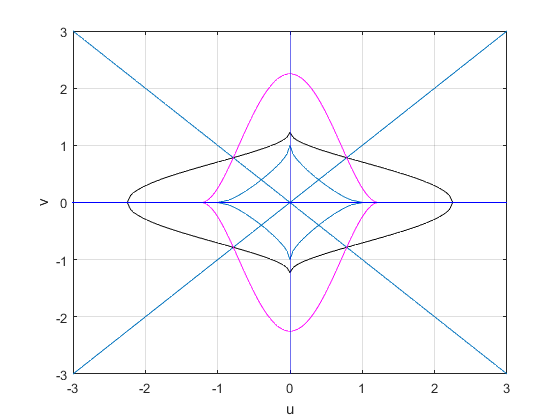}
		\caption{Parameters' domains for $p=1.5$, $\omega=1.5$}\label{w15p15}
	\end{center}
\end{figure}
\\
\begin{figure}[htbp]
	\begin{center}
		\includegraphics[scale=0.65]{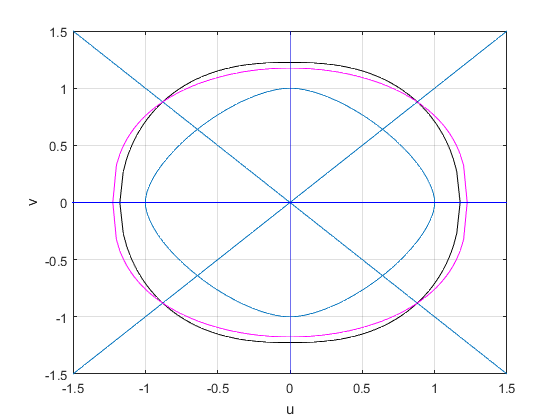}
		\caption{Parameters' domains for $p=3.5$, $\omega=1.5$}\label{w15p35}
	\end{center}
\end{figure}
\\

\end{document}